\documentclass[a4paper,12pt,reqno]{amsart}
\usepackage{mathtools}
\usepackage{amsfonts,amsmath}
\usepackage{amssymb,amsthm}		
\usepackage{amsthm}
\usepackage[T1]{fontenc}
\usepackage{lmodern}
\usepackage{hyperref}
\usepackage[all]{hypcap}
\usepackage{fullpage}
\usepackage{color}
\usepackage{graphicx}
\usepackage[dvipsnames]{xcolor}
\usepackage{fancyhdr}
\usepackage{psfrag}

\usepackage{indentfirst}
\usepackage{float}
\usepackage{latexsym}
\usepackage{graphics}
\usepackage{verbatim}

\newtheorem{theorem}{Theorem}
\newtheorem*{theorem*}{Theorem}
\newtheorem{lemma}[theorem]{Lemma}
\newtheorem{corollary}{Corollary}

\newtheorem{proposition}{Proposition}
\theoremstyle{definition}

\theoremstyle{remark}

\begin{document}

\author{Stephan Baier, Arpit Bansal}
\address{Stephan Baier, Ramakrishna Mission Vivekananda Educational Research Institute, Department of Mathematics, 
G. T. Road, PO Belur Math, Howrah, West Bengal 711202, 
India; email: email$_{-}$baier@yahoo.de}
\address{Arpit Bansal, Jawaharlal Nehru University, School of Physical Sciences, New-Delhi 110067, India; email: apabansal@gmail.com}
\title{Large sieve with sparse sets of moduli for $\mathbb{Z}[i]$}
\begin{abstract} We establish a general large sieve inequality with sparse sets $\mathcal{S}$ of modu\-li in the Gaussian integers
which are in a sense well-distributed in arithmetic progressions. This extends earlier work of S. Baier on the large sieve with sparse sets of moduli. 
We then use this result to obtain large sieve bounds for the cases when $\mathcal{S}$ consists of squares of Gaussian integers and of Gaussian primes. 
Our bound for the case of square moduli improves our recent result in \cite{BB}.
\end{abstract}
\maketitle

\section{Introduction}
The classical large sieve inequality with additive characters asserts that
\begin{equation*} \label{Classical}
\sum\limits_{q\le Q} \sum\limits_{\substack{a=1\\ (a,q)=1}}^q \left| \sum\limits_{M<n\le M+N} a_n e\left(n\cdot \frac{a}{q}\right) \right|^2
\le \left(Q^2+N-1\right)\sum\limits_{M<n\le M+N} |a_n|^2,
\end{equation*}
where $Q$, $N$ $\in \mathbb{N}$, $M \in \mathbb{Z}$ and $\{a_n\}$ is any arbitrary sequence of complex numbers. There are numerous applications of this inequality in analytic number theory, 
in particular, in sieve theory and to questions regarding the distribution of arithmetic functions in arithmetic progressions.\\ \\
The large sieve with sparse sets of moduli $q$, in particular with prime moduli and with square moduli were investigated by Wolke, Zhao and the first-named author in a series of papers 
(see \cite{Bai}, \cite{Ba2}, \cite{BZ} and \cite{Wol}). In the case of prime moduli, it was 
established by D. Wolke \cite{Wol} that
\begin{equation}\label{primemodZ}
\sum\limits_{p\le Q}\sum\limits_{a=1}^{p-1}\left|\sum\limits_{M<n\le M+N}a_ne\left(n\cdot \frac{a}{p}\right)\right|^2\le\frac{C}{1-\delta}\frac{Q^2\log\log Q}{\log Q}\sum\limits_{M<n\le M+N}|a_n|^2
\end{equation}  
provided that $ Q\geq 10$, $N = Q^{1+\delta}$, $0<\delta<1.$ Here, $C$ is an absolute constant.\\  
In the case of square moduli, it was first established by Zhao \cite{Zh1} that
\begin{equation}\label{squaremodZ1}
\begin{split}
& \sum\limits_{q\le Q}\sum\limits_{\substack{a = 1\\ (a,q) = 1}}^{q^2}\left|\sum\limits_{M<n\le M+N}a_ne\left(n\cdot\frac{a}{q^2}\right)\right|^2\\  \ll_\varepsilon&  (QN)^\varepsilon\left(Q^3+Q^2\sqrt{N}+\sqrt{Q}N\right)
\sum\limits_{M<n\le M+N}|a_n|^2.
\end{split}
\end{equation}
S. Baier improved this to
\begin{equation}\label{squaremodZ2}
\sum\limits_{q\le Q}\sum\limits_{\substack{a = 1\\ (a,q) = 1}}^{q^2}\left|\sum\limits_{M<n\le M+N}a_ne\left(n\cdot\frac{a}{q^2}\right)\right|^2 \ll_\varepsilon (QN)^\varepsilon\left(Q^3+Q^2\sqrt{N}+N\right)
\sum\limits_{M<n\le M+N}|a_n|^2.
\end{equation}
A further improvement was obtained by S.Baier and L. Zhao in \cite{BZ}, where the term $N+ Q^2\sqrt{N}$ was replaced by 
$N+\min\{\sqrt{Q}N,Q^2\sqrt{N}\}$. To date, this is the best known bound.
A generalization of the large sieve for number fields was established by M.Huxley \cite{Huxl}. For the number field $\mathbb{Q}(i)$, it takes the form
\begin{equation} \label{Huxley}
\sum\limits_{\substack{q\in \mathbb{Z}[i]\setminus\{0\}\\ \mathcal{N}(q)\le Q}} 
\sum\limits_{\substack{a \bmod{q}\\ (a,q)=1}} \left|\sum\limits_{\substack{n\in \mathbb{Z}[i]\\ \mathcal{N}(n)\le N}}  
a_n \cdot e\left(\Re\left(\frac{na}{q}\right)\right)\right|^2 \ll \left(Q^2+N\right)\sum\limits_{\substack{n\in \mathbb{Z}[i]\\ 
\mathcal{N}(n)\le N}} |a_n|^2.
\end{equation}
Here as in the following, $\mathcal{N}(q)$ denotes the norm of $q \in \mathbb{Z}[i]$, given by 
$$
        \mathcal{N}(q) := \Re (q)^2 + \Im (q)^2.
$$
In \cite{BB}, we studied the large sieve with square moduli for the number field $Q(i)$, i.e. we investigated the order of magnitude of the expression
     $$
     T := \sum\limits_{\substack{q\in \mathbb{Z}[i]\setminus\{0\}\\ \mathcal{N}(q)\le Q}} 
\sum\limits_{\substack{a \bmod{q^2}\\ (a,q)=1}} \left|\sum\limits_{\substack{n\in\mathbb{Z}[i]\\ \mathcal{N}(n)\le N}}a_n\cdot e\left(\Re\left(\frac{na}{q^2}\right)\right)\right|^2.
$$
We established an analogue of \eqref{squaremodZ1}, namely the inequality
\begin{equation} \label{earlier}
   T\ll_\varepsilon (QN)^\varepsilon \left(Q^3 + Q^2\sqrt{N} + \sqrt{Q}N\right)Z.
\end{equation} 
For comparison, Huxley's version \eqref{Huxley} of the large sieve in $\mathbb{Z}[i]$ with the set of moduli extended to all $q$ with $0<\mathcal{N}(q)\le Q$ implies only the bound
\begin{equation} \label{Huxleyapp}
           T\ll \left(N + Q^4\right)Z,
\end{equation}
which is weaker than \eqref{earlier} if $Q\gg N^{2/7+\varepsilon}$. On the other hand, it is easy to show that
\begin{equation} \label{individual}
\sum\limits_{\substack{a \bmod{q^2}\\ (a,q) = 1}} \left|\sum\limits_{\substack{n\in\mathbb{Z}[i]\\ \mathcal{N}(n)\le N}}a_n\cdot e\left(\Re\left(\frac{na}{q^2}\right)\right)\right|^2 \ll \left(N + \mathcal{N}(q)^2\right)Z,
\end{equation}
(in particular, this follows from our later Theorem \ref{lsR2}  with $\Delta = 1/\mathcal{N}(q)^2$ and $(x_r)$ being the sequence formed by all Farey fractions $a/q^2$ with 
$1\le |a|\le |q|^2$ and $(a,q) = 1$), which implies the bound 
\begin{equation} \label{summing}
T \ll Q\left(N + Q^2\right)Z
\end{equation}
 by summing up \eqref{individual} over all $q\in \mathbb{Z}[i]$ with $\mathcal{N}(q)\le Q$. This bound is weaker than \eqref{earlier} if $Q\ll N^{1/2-\varepsilon}$. Thus, 
 \eqref{earlier} is sharper than both \eqref{Huxleyapp} and \eqref{summing} if $N^{2/7+\varepsilon}\ll Q\ll N^{1/2-\varepsilon}$.
 
 The goal of this paper is to improve \eqref{Huxley} for sparse sets $\mathcal{S}$ of moduli which are in a sense well-distributed in arithmetic progressions. 
 As a consequence, we derive an analogue of \eqref{squaremodZ2} for 
 $\mathbb{Q}[i]$, thus improving \eqref{earlier}. Also, we establish a large sieve inequality with Gaussian prime moduli which is an analogue of \eqref{primemodZ}. Similarly as in \cite{BB},
 our method starts with an application of the large sieve for 
 $\mathbb{R}^2$. Then we convert the resulting counting problem back into one for $\mathbb{Q}(i)$. At this stage, we deviate significantly from the method in \cite{BB}, where we used
 Fourier analytic tools to attack the said counting problem. Instead, we proceed along similar lines as in \cite{Bai}, only using Diophantine approximation and elementary counting arguments.

\section{Main results}
Throughout this paper, we reserve the symbols $c_i$ $(i = 1,2,.....)$ for absolute constants and the symbol $\varepsilon$ for an arbitrary (small) positive number.
The $\ll$- constants in our estimates may depend on $\varepsilon$. As usual in analytic number theory, $\varepsilon$ may be different from line to line. We further 
suppose $(a_n)_{n \in \mathbb{Z}[i]}$ to be any sequence of complex numbers and $Q, N\in \mathbb{N}$. For $\alpha\in\mathbb{Z}[i]$, we set 
\begin{equation} \label{TrigPoly}
\begin{split}
     S(\alpha) := \sum\limits_{\substack{n \in \mathbb{Z}[i]\\ \mathcal{N}(n)\le N}} a_n\cdot e\left(\Re \left(n \alpha \right)\right)   
\end{split}
\end{equation} 
and 
\begin{equation}\label{Zdef}
 Z := \sum\limits_{\substack{n\in \mathbb{Z}[i]\\ \mathcal{N}(n)\le N}} |a_n|^2.
 \end{equation}
 
We further suppose that $\mathcal{S}$ $\subseteq$ $B(0,Q^{1/2}) \cap (\mathbb{Z}[i]\setminus \{0\})$, where $B(0,Q^{1/2})$ denotes the closed ball with center $0$ and radius $Q^{1/2}$, i.e.
$$
B(0,Q^{1/2}):=\{z\in \mathbb{C}\ :\ |z|\le Q^{1/2}\}.
$$
For $t\in \mathbb{Z}[i]\setminus\{0\}$ we put
 $$
 \mathcal{S}_t := \{q\in \mathbb{Z}[i] : tq \in \mathcal{S}\}.
 $$  
 We note that 
 $$
 \mathcal{S}_t \subset B\left(0,\frac{Q^{1/2}}{|t|}\right).
 $$
 We shall require that the number of elements of $\mathcal{S}_t$ in small regions of arithmetic progressions in $\mathbb{Z}[i]$ (which form shifted lattices in $\mathbb{C}$)
 does not differ too much from the expected 
 number. To measure the distribution of $\mathcal{S}_t$ in regions of arithmetic progressions, we define the quantity
 \begin{equation} \label{Atdef}
 A_t(u,k,l) := \sup\limits_{\substack{y\in \mathbb{C}\\  |y| \le \frac{\sqrt{Q}}{|t|}}}\left|\{ q\in \mathcal{S}_t\cap B(y,u) : q\equiv l\bmod{k}\}\right|,
 \end{equation}
 where $0\le u\le \sqrt{Q}/t$, $k\in \mathbb{Z}[i]\setminus\{0\}$ and $l\in \mathbb{Z}[i]$ with $(k,l) = 1$. Here $B(y,u)$ denotes the closed ball with center $y$ and radius $u$, i.e.
 $$
 B(y,u):=\{z\in \mathbb{C} : |z-y|\le u\}. 
 $$
 
We first establish the following large sieve inequality for general sets $\mathcal{S}$ of moduli in $\mathbb{Z}[i]$.
\begin{theorem} \label{result1} We have
\begin{equation*}
\begin{split}
 & \sum\limits_{q\in \mathcal{S}}  
 \sum\limits_{\substack{a \bmod{q}\\ (a,q)=1}} \left|S\left(\frac{a}{q}\right)\right|^2 \le \\ &  c_1 NZ \left( 1 + \sup\limits_{\substack{r\in \mathbb{Z}[i]\setminus\{0\}\\ 1 \le |r| \le N^{1/4}}} 
\sup\limits_{\substack{z\in \mathbb{C}\\ \frac{1}{N^{1/2}}\le |z| \le \frac{2}{|r|N^{1/4}}}} \sup\limits_{\substack{h \in \mathbb{Z}[i] \\ (h,r) = 1}} \sum\limits_{t|r}  
 \sum\limits_{\substack{0<|m|\le \frac{3|rz|\sqrt{Q}}{|t|} \\ (m,\frac{r}{t}) = 1}}
 A_t \left( \frac{\sqrt{Q}}{\sqrt{N}|zt|},\frac{r}{t},hm \right)\right),
\end{split}
\end{equation*}
where $Z$ is defined as in \eqref{Zdef}.
\end{theorem}

If we assume the set $\mathcal{S}_t$ to be nearly evenly distributed in the residue classes $l\bmod{k}$, then, if 
$B(y,u)\subset B(0,\sqrt{Q}/|t|)$, the expected cardinality of the set 
$$
\{q\in \mathcal{S}_t \cap B(y,u) : q\equiv l \bmod{k}\}
$$
is 
$$
        \asymp     \frac{|\mathcal{S}_t|/\mathcal{N}(k)}{Q/|t|^2}\cdot u^2.
$$
This suggests to set a condition of the form
\begin{equation} \label{condi}
     A_t(u,k,l) \le \left(1 + \frac{|\mathcal{S}_t|/\mathcal{N}(k)}{Q/|t|^2}\cdot u^2 \right)X,
\end{equation}
where $X\geq 1$ is thought to be small compared to $Q$ and $N$.
Under the condition \eqref{condi}, we shall infer the following bound from Theorem \ref{result1}.

\begin{theorem} \label{result2} 
Suppose the condition (11) to hold for all $t$, $k$, $l$, $u$ with $|t|\le N^{1/4}$, $|k|\le N^{1/4}/|t|$, $(k,l) = 1$ and 
$|k|\sqrt{Q}/(2N^{1/4}) \le u \le \sqrt{Q}/|t|$. Then 
\begin{equation} \label{gensparse}
\sum\limits_{q\in \mathcal{S}} \sum\limits_{\substack{a\bmod{q} \\ (a,q) = 1}} \left|S\left(\frac{a}{q}\right)\right|^2 \le 
c_2\left(N + QXN^{\varepsilon}\left(\sqrt{N} + |\mathcal{S}|\right)\right)Z.
\end{equation}
\end{theorem}

Inequality \eqref{gensparse} is stronger than the ``trivial bound''
\begin{equation} \label{triv}
\sum\limits_{q\in \mathcal{S}} \sum\limits_{\substack{a\bmod{q} \\ (a,q) = 1}} 
\left|S\left(\frac{a}{q}\right)\right|^2 \ll \left(Q^2+N\right)Z
\end{equation}
following directly from Huxley's large sieve \eqref{Huxley}, if 
$$
N^{1+\varepsilon} \ll QX\left(\sqrt{N} + |\mathcal{S}|\right)\ll Q^{2-\varepsilon}.
$$
Employing Theorem \ref{result2} with $\mathcal{S}$ a set of non-zero squares of norm $\le Q^2$, we shall derive the following improvement of \eqref{earlier}. 

\begin{theorem}\label{squaremoduli Zi} 
We have
 \begin{equation} \label{squaremods}
\sum\limits_{\substack{q\in \mathbb{Z}[i]\setminus\{0\}\\ \mathcal{N}(q)\le Q}} 
\sum\limits_{\substack{a \bmod{q^2}\\ (a,q)=1}} \left|\sum\limits_{\substack{n\in\mathbb{Z}[i]\\ \mathcal{N}(n)\le N}}a_n\cdot e\left(\Re\left(\frac{na}{q^2}\right)\right)\right|^2 
 \ll (QN)^{\varepsilon}\left(Q^3 + Q^2\sqrt{N} + N\right)Z,
 \end{equation}
 where $\varepsilon$ is any positive constant, and the implied constant $\ll$-constant depends only on $\varepsilon$.
 \end{theorem}
 
 This bound  is stronger than the three bounds \eqref{earlier}, \eqref{Huxleyapp} and \eqref{summing} if 
 $N^{1/4+\varepsilon}\ll Q\ll N^{1/3 - \varepsilon}$. 
 When $\mathcal{S}$ is the full set of all Gaussian primes with norm $\le Q$, we shall establish the following version of the large sieve for 
 $\mathbb{Z}[i]$.
 
 \begin{theorem}\label{primemoduli}
 Let $Q\geq 16$, $N = Q^{1+\delta}/16$, $0 < \delta < 1$. Then there is an absolute constant $c_1$ such that
 \begin{equation}
 \sum\limits_{ 
 \mathcal{N}(p)\le Q}\sum\limits_{\substack{a\bmod{p}\\ (a,p)=1}}\left|S\left(\frac{a}{p}\right)\right|^2 \le 
 \frac{c_1}{1-\delta}\frac{Q^2\log\log Q}{\log Q}Z,
 \end{equation}
 where $p$ runs over the Gaussian primes. 
 \end{theorem}
 
 \section{Large sieve for $\mathbb{R}^2$}
 We shall employ the following version of the large sieve for $\mathbb{R}^2$, proved below.
 \begin{theorem}\label{lsR2}
 Let $R, N\in \mathbb{N}$, $N\geq 2$, $x_1,x_2, ....,x_R \in \mathbb{R}^2$ and $(a_n)_{n\in \mathbb{Z}^2}$ be any double sequence of complex numbers. Suppose that $0<\Delta \le 1/2$. Put
 \begin{equation}
   K(\Delta) := \sup\limits_{\alpha\in \mathbb{R}^2} \left|\left\{r\in\{1,2,...,R\} : \min\limits_{z\in\mathbb{Z}^2}||x_r-\alpha-z||_2\le \Delta^{1/2}\right\}\right|,
\end{equation}
where $||.||_2$ is the Euclidean norm on $\mathbb{R}^2$.
Then,
\begin{equation}
\sum\limits_{r=1}^{R}\left|\sum\limits_{\substack{n\in \mathbb{Z}^2\\ ||n||_2^2\le N}}a_n\cdot e(n\cdot x_r)\right|^2 \le c_4K(\Delta)(N+ \Delta^{-1})Z.
\end{equation}   
\end{theorem}

Here as in the following, $||x||_2$ denotes the Euclidean norm of $s\in \mathbb{R}^2$, given by
$$       
||(x_1,x_2)||_2 = \sqrt{x_1^2+x_2^2}.
$$
To prove Theorem 5, we use the duality principle and the Poisson summation formula for $\mathbb{R}^2$.
\begin{proposition}[Duality principle, Theorem 288 in \cite{HLP}] \label{duality} Let $C = [c_{mn}]$ be a finite matrix with complex entries. 
The following two statements are equivalent: 
\begin{enumerate}
\item For any complex numbers $a_n$, we have 
\begin{align*}
\sum_m \mathrel \Big |\sum_n a_n c_{mn}\Big |^2 \leq \Delta \sum_n |a_n|^2.
\end{align*}
\item For any complex numbers $b_m$, we have
\begin{align*}
\sum_n \mathrel \Big |\sum_m b_m c_{mn}\Big |^2 \leq \Delta \sum_m |b_m|^2. 
\end{align*}
\end{enumerate}
\end{proposition}

\begin{proposition}[Poisson summation formula, see \cite{StW}] \label{poisson} 
Let $f:\mathbb{R}^2 \rightarrow \mathbb{C}$ be a smooth function of rapid decay and $\Lambda$ be a lattice 
of full rank in $\mathbb{R}^2$. Then
$$
\sum\limits_{y\in \Lambda} f(y) = \frac{1}{\mbox{\rm Vol}(\mathbb{R}^2/\Lambda)} \cdot \sum\limits_{x\in \Lambda'} \hat{f}(x),
$$
where $\Lambda'$ is the dual lattice, $\hat{f}$ is the Fourier transform of $f$, defined as 
$$
\hat{f}(x)=\int\limits_{\mathbb{R}^2} f(y)e\left(-x\cdot y\right)\ dy,
$$
and $\mbox{\rm Vol}(\mathbb{R}^2/\Lambda)$ is the volume of a fundamental mesh of $\Lambda$.
\end{proposition}

Here as in the following, by {\it rapid decay} we mean that the function $f:\mathbb{R}^2 \rightarrow \mathbb{C}$ in question satisfies the bound
$$
f(y)\ll \left(1+||y||_2\right)^{-C}
$$
for some $C>2$.\\ \\
{\bf Proof of Theorem \ref{lsR2}:}
Note that
\begin{equation} \label{note}
\begin{split}
K(\Delta)  = & \sup\limits_{\alpha\in\mathbb{R}^2}\left|\left\{r\in\{1,2,...,R\} : \min\limits_{z\in\mathbb{Z}^2}||x_r-\alpha-z||_2\le\Delta^{1/2}\right\}\right| 
\\ \ge & \sup\limits_{\alpha\in\mathbb{R}^2}\left|\left\{r\in\{1,2,...,R\} : \max\{||x_r^{(1)}-\alpha^{(1)}||, ||x_r^{(2)}-\alpha^{(2)}||\}\le 
\frac{\Delta^{1/2}}{\sqrt{2}}\right\}\right|,  
\end{split}
\end{equation}
where $||u||$ is the distance of $u\in \mathbb{R}$ to the nearest integer and we write
$$  x_r = \left(x_r^{(1)},x_r^{(2)}\right), \alpha = \left(\alpha^{(1)},\alpha^{(2)}\right) \ \text{and} \ z = \left(z^{(1)},z^{(2)}\right)$$
for $r = 1,2,...,R$. 
Now, let $S = \{x_1,x_2,...,x_R\}$. Taking Proposition \ref{duality}, the duality principle, into account, it suffices to prove that
\begin{equation*}
\sum\limits_{\substack{n\in\mathbb{Z}^2\\||n||_2\le N^{1/2}}}\left|\sum\limits_{x\in S}b_x\cdot e(n\cdot x)\right|^2\le c_5K(\Delta)(N+ \Delta^{-1})\sum\limits_{x\in S}|b_x|^2
\end{equation*}
for any complex numbers $b_x$. To this end, for $x = \left(x^{(1)},x^{(2)}\right)\in \mathbb{R}^2$, we define 
$$
\phi(x) = \left(\frac{\sin(\pi x^{(1)})}{2x^{(1)}}\right)^2\left(\frac{\sin(\pi x^{(2)})}{2x^{(2)}}\right)^2.
$$
It is clear that $\phi(x)\geq 0$ for all $x\in \mathbb{R}^2$ and $\phi(x)\geq 1$ if $|x^{(1)}|, |x^{(2)}|\le 1/2$. Moreover, the Fourier transformation of $\phi(x)$ equals
$$
\hat{\phi}(s) = \left(\frac{\pi^2}{4}\right)^2\max\left\{1-|s^{(1)}|,0\right\}\cdot\max\left\{1-|s^{(2)}|,0\right\},
$$
where $s = \left(s^{(1)},s^{(2)}\right)\in \mathbb{R}^2$. Hence,
\begin{equation*}
\begin{split}
\sum\limits_{\substack{n\in\mathbb{Z}^2\\||n||_2\le N^{1/2}}}\left|\sum\limits_{x\in S}b_x\cdot e(n\cdot x)\right|^2 & \le 
\sum\limits_{\substack{n\in \mathbb{Z}^2\\ ||n||_2< N^{1/2} + \Delta^{-1/2}}}\left|\sum\limits_{x\in S}b_xe(n\cdot x)\right|^2 \\ &
\le \sum\limits_{n\in \mathbb{Z}^2} \phi\left(\frac{n}{2(N^{1/2}+\Delta^{-1/2})}\right)\left|\sum\limits_{x\in S}b_xe(n\cdot x)\right|^2 \\ &
 =  \sum\limits_{n\in \mathbb{Z}^2} \phi\left(\frac{n}{2(N^{1/2}+\Delta^{-1/2})}\right)\sum\limits_{x,x'\in S}b_x\overline{b_{x'}}e\left(n(x-x')\right) \\ &
 =  \sum\limits_{x,x'\in S} b_x\overline{b_{x'}}V(x-x'),
\end{split}
\end{equation*}
where
$$ V(y) = \sum\limits_{n\in \mathbb{Z}^2} \phi\left(\frac{n}{2(N^{1/2}+\Delta^{-1/2})}\right)e(n\cdot y).$$
Using Proposition \ref{poisson}, the Poisson summation formula, we transform $V(y)$ into
\begin{equation*}
\begin{split}
V(y) = & 4(N^{1/2} + \Delta^{-1/2})^2\sum\limits_{\alpha\in {y+\mathbb{Z}^2}}\tilde{\phi}\left(2(N^{1/2} + \Delta^{-1/2})\alpha\right) \\ 
       = & 4(N^{1/2} + \Delta^{-1/2})^2\sum\limits_{\alpha\in {-y+\mathbb{Z}^2}}\hat{\phi}\left(2(N^{1/2} + \Delta^{-1/2})\alpha\right)  \\ 
       = & 4(N^{1/2} + \Delta^{-1/2})^2\left(\frac{\pi^2}{4}\right)^2\max\left\{1-|2(N^{1/2}+\Delta^{-1/2})y^{(1)}|,0\right\}\cdot \\ 
       & \max\left\{1-|2(N^{(1/2)} + \Delta^{-1/2})y^{(2)}|,0\right\},
 \end{split}
\end{equation*}
where $\tilde{\phi}$ is inverse Fourier transform and $\hat{\phi}$ is the Fourier transform of $\phi$. Therefore,
\begin{equation*}
\begin{split}
\sum\limits_{\substack{n\in\mathbb{Z}^2\\||n||_2\le N^{1/2}}}\left|\sum\limits_{x\in S}b_x\cdot e(n\cdot x)\right|^2 \le &\frac{\pi^4}{16}\left(N^{1/2} + \Delta^{-1/2}\right)^2 
\sum\limits_{\substack{x,x'\in S\\ ||x^{(1)}-x'^{(1)}||\le \frac{1}{2\left(N^{1/2}+\Delta^{-1/2}\right)}
\\ ||x^{(2)}-x'^{(2)}||\le \frac{1}{2\left(N^{1/2}+\Delta^{-1/2}\right)}}} |b_x||b_x'|\\
\le &\frac{\pi^4}{16}\left(N^{1/2} + \Delta^{-1/2}\right)^2 
\sum\limits_{\substack{x,x'\in S\\ ||x^{(1)}-x'^{(1)}||\le \frac{\Delta^{1/2}}{2}
\\ ||x^{(2)}-x'^{(2)}||\le \frac{\Delta^{1/2}}{2}}} |b_x||b_x'|.
\end{split}
\end{equation*}
Now we observe that
      $$ |b_x||b_x'|\le \frac{1}{2}\cdot\left(|b_x|^2+|b_x'|^2\right).$$
Using \eqref{note}, it follows that
      $$\sum\limits_{\substack{n\in\mathbb{Z}^2\\||n||_2\le N^{1/2}}}\left|\sum\limits_{x\in S}b_x\cdot e(n\cdot x)\right|^2
      \le c_4K(\Delta)\left(N+\Delta^{-1}\right)Z.$$
This completes the proof. $\Box$

\section{Conversion into a counting problem}
Now we return to the large sieve for $\mathbb{Q}(i)$. We aim to estimate the quantity
$$
U:=\sum\limits_{q\in \mathcal{S}} 
\sum\limits_{\substack{a \bmod{q}\\ (a,q)=1}} \left|\sum\limits_{\substack{n\in \mathbb{Z}[i]\\ \mathcal{N}(n)\le N}}  
a_n \cdot e\left(\Re\left(\frac{na}{q}\right)\right)\right|^2. 
$$
Our first step is to re-write $U$ in the form
\begin{equation} \label{U}
U=\sum\limits_{q\in \mathcal{S}} \sum\limits_{\substack{a \bmod{q}\\ (a,q)=1}} \left|\sum\limits_{\substack{n\in \mathbb{Z}[i]\\ \mathcal{N}(n)\le N}}  
a_n \cdot e\left(\left(\frac{xu+yv}{\mathcal{N}(q)},\frac{xv-yu}{\mathcal{N}(q)}\right)\cdot (s,t)\right)\right|^2, 
\end{equation}
where
$$
q=u+iv,\quad a=x+iy,\quad n=s+it. 
$$
To bound $U$, we employ Theorem \ref{lsR2}, which immediately gives us the following.

\begin{corollary} \label{preform} For $U$ as defined in \eqref{U}, we have the bound
$$
U\ll K(\Delta)(N+\Delta^{-1})Z,
$$
where $Z$ is defined as in \eqref{Zdef} and 
\begin{equation*}
\begin{split}
 K(\Delta):= \sup\limits_{\alpha\in\mathbb{R}^2}
 \Bigg|\Bigg\{ & (a,q)\in \mathbb{Z}[i]\times \mathcal{S} : (a,q)=1,\\ & \min\limits_{z\in \mathbb{Z}^2} 
 \left|\left|\left(\frac{xu+yv}{\mathcal{N}(q)},\frac{xv-yu}{\mathcal{N}(q)}\right)-\alpha -z\right|\right|_2\le \Delta^{1/2}\Bigg\} \Bigg|.
 \end{split}
\end{equation*}
\end{corollary}
Thus, we have converted the problem into a counting problem in $\mathbb{R}^2$, 
which we shall now interpret as a counting problem in $\mathbb{C}$. We observe that
$$
\frac{\overline{a}}{\overline{q}}=\frac{x-iy}{u-iv}=\frac{xu+yv}{\mathcal{N}(q)}+\frac{xv-yu}{\mathcal{N}(q)}i.
$$
It follows that
\begin{equation*}
\begin{split}
K(\Delta) = & \sup\limits_{\alpha\in\mathbb{C}}\left|\left\{(a,q) \in \mathbb{Z}[i]\times \mathcal{S} :  (a,q) = 1,\ 
\min\limits_{z\in \mathbb{Z}[i]} \left|\frac{\overline{a}}{\overline{q}}-\alpha-z \right|\le \Delta^{1/2}\right\}\right| \\
= & \sup\limits_{\alpha\in\mathbb{C}} P(\alpha),
\end{split}
\end{equation*}
where  
\begin{equation}
\begin{split}
P(\alpha) := & \left|\left\{(a,q) \in \mathbb{Z}[i]\times \mathcal{S} :  (a,q) = 1,\ 
\left| \frac{a}{q}-\alpha\right|\le \Delta^{1/2}\right\}\right|\\
           = & \sum\limits_{\substack{q\in \mathcal{S}, (a,q) = 1\\ \frac{a}{q}\in B(\alpha,\Delta^{1/2})}} 1. 
\end{split}
  \end{equation}      
Now we are left with counting Farey fractions in $\mathbb{C}$.
  
 \section{Counting Farey fractions in small regions in $\mathbb{C}$} 
To estimate $P(\alpha)$, we approximate $\alpha$ by a suitable element of $\mathbb{Q}(i)$. Let
\begin{equation}
            \tau := \frac{1}{\Delta^{1/4}}.
\end{equation} 
Then, using the Dirichlet approximation theorem in $\mathbb{C}$ (see \cite{MMD}), $\alpha$ can be written in the form
\begin{equation} \label{approxi}
\alpha = \frac{b}{r} + z, \ \text{where} \  b,r\in \mathbb{Z}[i],\ (b,r) = 1,\ |z|<\frac{2}{|r|\tau},\ 0<\ |r|\le\tau.
\end{equation}
Thus, it suffices to estimate $P(b/r+z)$ for all $b$, $r$, $z$ satisfying \eqref{approxi}.
We further note that we can restrict ourselves to the case when
\begin{equation} \label{zcond}
|z|\geq \Delta^{1/2}
\end{equation}
because if $|z|<\Delta^{1/2}$, then
$$ P(\alpha)\le P\left(\frac{b}{r} + \Delta^{1/2}\right) + P\left(\frac{b}{r} - \Delta^{1/2}\right) + P\left(\frac{b}{r} +i\Delta^{1/2}\right) + P\left(\frac{b}{r}-i\Delta^{1/2}\right).$$
We deduce the following.

\begin{lemma} \label{cr1}
We have
\begin{equation}
K(\Delta) \le 4\sup\limits_{\substack{r\in\mathbb{Z}[i]\\ 1\le |r|\le \tau}}\sup\limits_{\substack{b\in\mathbb{Z}[i]\\ (b,r) =1}}\sup\limits_{\substack{z\in\mathbb{C}\\ 
\Delta^{1/2}\le |z|\le \frac{2}{|r|\tau}}}P\left(\frac{b}{r}+z\right).
\end{equation}
\end{lemma}

The next lemma provides a first estimate for $P\left(\frac{b}{r}+z\right)$.
\begin{lemma}\label{cr2}
Suppose that the conditions \eqref{approxi} and \eqref{zcond} are satisfied. Suppose further that 
\begin{equation}
\frac{\sqrt{Q}\Delta^{1/2}}{|z|}\le \delta^{1/2}\le \sqrt{Q}.
\end{equation}
Let 
\begin{equation} \label{Jdef}
J(\mu,\nu) := B(rz\mu,2|rz|\nu)
\end{equation}
and
\begin{equation} \label{Pidef}
\Pi (y,\delta) := \sum\limits_{q\in\mathcal{S}\cap B(y,\delta^{1/2})}\sum\limits_{\substack{m\in J(y,\delta^{1/2})\\ m\equiv -bq\bmod{r} \\ m\neq 0}} 1.
\end{equation}
Then,
$$
P\left(\frac{b}{r} + z\right) \le 16 + \frac{4}{\pi \delta}\int\limits_{B(0,\sqrt{Q})}\Pi (y,\delta) \ dy,
$$
where $B(0,\sqrt{Q})$ is the closed ball with center $0$ and radius $\sqrt{Q}$.
\end{lemma}

{\bf Proof of Lemma \ref{cr2}:} 
Define
\begin{equation*}
P(\alpha, y, \delta) := \sum\limits_{\substack{q\in\mathcal{S}\cap B(y,\delta^{1/2})\\ (a,q) =1\\  \frac{a}{q}\in B(\alpha,\Delta^{1/2})}} 1.
\end{equation*}
Then, if $\delta\le Q$, we have
\begin{equation*}
\begin{split}
\int\limits_{B(0,\sqrt{Q})}P(\alpha, y, \delta) \ dy 
= & \sum\limits_{\substack{q\in\mathcal{S} \\ (a,q) = 1\\ \frac{a}{q}\in B(\alpha,\Delta^{1/2})}}
\int\limits_{B(0,\sqrt{Q})\cap B(q,\delta^{1/2})} dy \\ 
\geq & \frac{\pi \delta}{4}\sum\limits_{\substack{q\in\mathcal{S}\\ (a,q) =1\\ \frac{a}{q}\in B(\alpha,\Delta^{1/2})}}1 = \frac{\pi \delta}{4} P(\alpha).
\end{split}
\end{equation*} 
This implies
\begin{equation} \label{Palpha}
 P(\alpha) \le \frac{4}{\pi \delta}\int\limits_{D(0,\sqrt{Q})}P(\alpha,y,\delta)\ dy
\end{equation}
whenever $\delta \le Q$.

Now, for $a/q\in B(\alpha,\Delta^{1/2})$, we have
$$
|a-q\alpha|\le |q|\Delta^{1/2}.
$$
From this and $\alpha = b/r+z$, we deduce that
$$
  \left|a-q\frac{b}{r}-qz\right|\le |q|\Delta^{1/2}
$$
and hence
\begin{equation*}
  |ar-bq-qrz| \le \sqrt{Q}\Delta^{1/2}|rz|/|z|.
\end{equation*}
If $\sqrt{Q}\Delta^{1/2}/|z|\le \delta^{1/2}$, then it follows that
\begin{equation}\label{inequality1}
|ar-bq-qrz| \le \delta^{1/2}|rz|,    
\end{equation}
and if $q\in B(y,\delta^{1/2})$, then 
\begin{equation}\label{inequality2}
 |qrz- yrz| \le \delta^{1/2}|rz|.
\end{equation}
From \eqref{inequality1} and \eqref{inequality2}, we deduce that
$$
 |ar-bq-yrz|\le 2\delta^{1/2}|rz|. 
$$
If $ar-bq = 0$, then $q$ is associated to $r$ (we write $q \approx r$) because $(a,q) = 1 = (b,r)$. Writing $m=ar-bq$ and recalling \eqref{Jdef} and \eqref{Pidef}, we deduce that
\begin{equation*}
\begin{split}
P(\alpha,y,\delta) = & \sum\limits_{\substack{q\in\mathcal{S}\cap B(y,\delta^{1/2})\\ (a,q) =1\\ \frac{a}{q}\in B(\alpha,\Delta^{1/2})}}1 = 
\sum\limits_{q\in\mathcal{S}\cap B(y,\delta^{1/2})}\sum\limits_{\substack{(a,q)=1\\ \frac{a}{q}\in B(\alpha,\Delta^{1/2})}}1 \\ 
     \le & \sum\limits_{q\in\mathcal{S}\cap B(y,\delta^{1/2})}\sum\limits_{\substack{m\in J(y,\delta^{1/2})\\ m\equiv -bq\bmod{r}}}1 \\ 
     = & \sum\limits_{\substack{q\in\mathcal{S}\cap B(y,\delta^{1/2})\\ q\approx r}} 1 + \Pi(y,\delta).
\end{split}
\end{equation*}
Hence, from \eqref{Palpha}, we obtain
\begin{equation*}
\begin{split}
P(\alpha) & = P\left(\frac{b}{r}+z\right)\le \frac{4}{\pi \delta}\int\limits_{B(0,\sqrt{Q})}\sum\limits_{\substack{q\in\mathcal{S}\cap B(y,\delta^{1/2}) \\ q \approx r}} 1\ dy + 
\frac{4}{\pi \delta}\int\limits_{B(0,\sqrt{Q})}\Pi(y,\delta)\ dy \\ &
         \le \frac{4}{\pi \delta} \sum\limits_{q\approx r} \int\limits_{B(q,\delta^{1/2})}1\ dy + \frac{4}{\pi \delta}\int\limits_{B(0,\sqrt{Q})}\Pi(y,\delta)\ dy \\ & 
         = 16 +\frac{4}{\pi \delta}\int\limits_{B(0,\sqrt{Q})}\Pi(y,\delta)\ dy
\end{split}
\end{equation*}
since $r$ has precisely 4 associates.
This completes the proof. $\Box$

\section{Proofs of Theorems 1 and 2}
Next, we express $\Pi(y,\delta)$ in terms of $A_t(u,k,l)$. This shall lead us to the following estimate for $P(b/r+z)$.  

\begin{lemma}\label{cr3}
We have
\begin{equation} \label{PAt}
P\left(\frac{b}{r}+z\right)\le 16 + c_3\sum\limits_{t|r}\sum\limits_{\substack{0<|m|\le \frac{3|rz|\sqrt{Q}}{|t|}\\(m,\frac{r}{t})=1}} 
A_t\left(\frac{\sqrt{Q\Delta}}{|z||t|}, \frac{r}{t},-\overline{b}m\right),
\end{equation}
where $\overline{b}b\equiv 1\bmod{r}$.
\end{lemma}

On choosing $\Delta := 1/N$, Theorem 1 follows immediately from Corollary \ref{preform}, Lemma \ref{cr1} and Lemma \ref{cr3} .\\ \\
{\bf Proof of Lemma \ref{cr3}:}
We split $\Pi(y,\delta)$ into
\begin{equation*}
\Pi(y,\delta) = 
\sum\limits_{\substack{t\\(q,r) \approx t}}\sum\limits_{q\in\mathcal{S}\cap B(y,\delta^{1/2})}
\sum\limits_{\substack{m\in J(y,\delta^{1/2})\\ m\equiv -bq\bmod{r}\\ m\neq 0}}1,
\end{equation*} 
where $t$ runs over a maximal set of mutually non-associate elements of $\mathbb{Z}[i]\setminus \{0\}$, and $(q,r)\approx t$ means that $t$ is a greatest common divisor of $(q,r)$ (unique up to
associates). 
Writing $\tilde{q}:=q/t$ and $\tilde{m}:=m/t$, it follows that
\begin{equation*}
\begin{split}
\Pi(y,\delta)\le & 
\sum\limits_{t|r}\sum\limits_{\substack{\frac{q}{t}\in\mathcal{S}_t\cap B\left(\frac{y}{t},\frac{\delta^{1/2}}{|t|}\right)\\(\frac{q}{t},\frac{r}{t}) =1}}
\sum\limits_{\substack{\frac{m}{t}\in J\left(\frac{y}{t},\frac{\delta^{1/2}}{|t|}\right)\\ \frac{m}{t}\equiv -b\frac{q}{t}\bmod{\frac{r}{t}}\\ 
\frac{m}{t}\neq 0}} 1\\ 
= & 
\sum\limits_{t|r}\sum\limits_{\substack{\tilde{q}\in\mathcal{S}_t\cap B\left(\frac{y}{t},\frac{\delta^{1/2}}{|t|}\right)\\ 
(\tilde{q},\frac{r}{t})=1}} \sum\limits_{\substack{\tilde{m}\in J\left(\frac{y}{t},\frac{\delta^{1/2}}{|t|}\right)\\
             \tilde{m}\equiv -b\tilde{q}\bmod{\frac{r}{t}}\\ \tilde{m}\neq 0}}1 \\ 
= &
\sum\limits_{t|r}\sum\limits_{\substack{\tilde{m}\in J\left(\frac{y}{t},\frac{\delta^{1/2}}{|t|}\right)\\ 
(\tilde{m},\frac{r}{t})=1\\ \tilde{m}\neq 0}}\sum\limits_{\substack{\tilde{q}\in\mathcal{S}_t\cap B\left(\frac{y}{t},\frac{\delta^{1/2}}{|t|}\right)\\
             \tilde{q}\equiv -\overline{b}\tilde{m}\bmod{\frac{r}{t}}}} 1. 
\end{split}
\end{equation*}
Hence, by definition of $A_t(u,k,l)$ in \eqref{Atdef}, we have 
\begin{equation*}
\Pi(y,\delta)
\le \sum\limits_{t|r}\sum\limits_{\substack{m\in J\left(\frac{y}{t},\frac{\delta^{1/2}}{|t|}\right)\\ (m,\frac{r}{t})=1\\ m\neq  0}}
A_t\left(\frac{\delta^{1/2}}{|t|},\frac{r}{t},-\overline{b}m\right).
\end{equation*}
Integrating the last line over $y$ in the ball $B(0,\sqrt{Q})$ and rearranging the order of summation and integration, we obtain
\begin{equation*}
\begin{split}
\int\limits_{B(0,\sqrt{Q})}\Pi(y,\delta^{1/2})\ dy & \le \sum\limits_{t|r} \int\limits_{B(0,\sqrt{Q})}\sum\limits_{\substack{m\in 
J\left(\frac{y}{t},\frac{\delta^{1/2}}{|t|}\right)\\ (m,\frac{r}{t}) =1 \\ m\neq 0}}
     A_t\left(\frac{\delta^{1/2}}{|t|},\frac{r}{t},-\overline{b}m\right) \ dy \\ &
    \le \sum\limits_{t|r}\int\limits_{B(0,\sqrt{Q})}\sum\limits_{\substack{|m-\frac{yrz}{t}|\le 2\delta^{1/2}\frac{|rz|}{|t|}\\ (m,\frac{r}{t}) =1\\ m\neq 0}}
     A_t\left(\frac{\delta^{1/2}}{|t|},\frac{r}{t},-\overline{b}m\right)\ dy \\ & 
     \le c_4\delta\sum\limits_{t|r}\sum\limits_{\substack{0<|m|\le (\sqrt{Q}+2\delta^{1/2})\frac{|rz|}{|t|}\\ (m,\frac{r}{t})=1}} A_t\left(\frac{\delta^{1/2}}{|t|},\frac{r}{t},
     -\overline{b}m\right).
\end{split}
\end{equation*}
Choosing 
$$
\delta := \frac{Q\Delta}{|z|^2}
$$
and using Lemma \ref{cr2}, we obtain \eqref{PAt}. $\Box$\\ \\
{\bf Proof of Theorem \ref{result2}:}   
From equation \eqref{condi}, we get 
\begin{equation*}
\begin{split}
& \sum\limits_{t|r}\sum\limits_{\substack{0<|m|\le \frac{3|rz|\sqrt{Q}}{|t|}\\ (m,\frac{r}{t}) = 1}} 
A_t\left(\frac{\sqrt{Q}}{\sqrt{N}|zt|},\frac{r}{t},hm \right) \\
\le & \sum\limits_{t|r}\sum\limits_{\substack{0<|m|\le \frac{3|rz|\sqrt{Q}}{|t|}\\ (m,\frac{r}{t}) = 1}} 
\left(1 + \frac{|\mathcal{S}_t|/\mathcal{N}(r/t)}{Q/|t|^2}\cdot \frac{Q}{N|zt|^2}\right)X 
\\ \le & \sum\limits_{t|r}\left(1 + \frac{|\mathcal{S}_t|\mathcal{N}(t) |t|^2 Q}{\mathcal{N}(r) QN|zt|^2}\right)\cdot \frac{9|rz|^2Q}{|t|^2} \cdot X \\ 
= & 9\sum\limits_{t|r} \left(\frac{|rz|^2}{|t|^2}\cdot N + |\mathcal{S}_t|\right)\frac{QX}{N}. 
\end{split}
\end{equation*}  
We deduce that the right-hand side of the inequality in Theorem \ref{result1} is dominated by
\begin{equation*}
\begin{split}
& c_1NZ\left(1 + 9\sup\limits_{\substack{r\in\mathbb{Z}[i]\\ 1\le|r|\le N^{1/4}}}\sum\limits_{t|r}\left(\frac{4N}{\sqrt{N}|t|^2} + 
|\mathcal{S}|\right)\frac{QX}{N}\right)\\ 
\le & c_2 \left(N + QXN^\varepsilon \left(\sqrt{N} + |\mathcal{S}|\right)\right)Z.
 \end{split}
 \end{equation*}
This completes the proof. $\Box$

\section{Proof of Theorem \ref{squaremoduli Zi}}
In this section, we derive Theorem \ref{squaremoduli Zi} from Theorem \ref{result2}. First, we rewrite the sum in question in the form
$$
T = \sum\limits_{\substack{q\in \mathbb{Z}[i]\\ \mathcal{N}(q)\le Q}}\sum\limits_{\substack{a\bmod{q^2}\\ (a,q)=1}}\left|S\left(\frac{a}{q^2}\right)\right|^2  
   =\sum\limits_{q\in\mathcal{S}}\sum\limits_{\substack{a\bmod{q}\\ (a,q)=1}}\left|S\left(\frac{a}{q}\right)\right|^2,$$
where $\mathcal{S}$ is the set of non-zero squares with norm $\le Q^2$. 
We split up the set $\mathcal{S}$ into $\mathcal{O}(\log 2Q)$ subsets of the form
\begin{equation*}
\mathcal{S}(Q_0) = \{q\in\mathcal{S} : Q_0<\mathcal{N}(q)\le 2Q_0\},
\end{equation*}
where $2Q_0\le Q^2$. Then we shall apply Theorem \ref{result2} to bound the quantity
$$
\sum\limits_{q\in\mathcal{S}(Q_0)}\sum\limits_{\substack{a\bmod{q}\\ (a,q)=1}}\left|S\left(\frac{a}{q}\right)\right|^2.
$$

As previously, we define
\begin{equation*}
\mathcal{S}_t(Q_0) := \{ q\in \mathbb{Z}[i] : tq\in \mathcal{S}(Q_0)\}.
\end{equation*}
Let $t = \epsilon p_1^{v_1} p_2^{v_2}... p_n^{v_n}$ be a prime factorization of $t$ in $\mathbb{Z}[i]$ (unique up to associates of $p_1,...,p_n$ and the unit $\epsilon$). For $i = 1,2,...,n$ let
$$
u_i = 
      \begin{cases}
      v_i, & \text{if} \  \ 2|{v_i} \\
      v_i+1, &\text{if} \ \ 2\nmid {v_i}.
      \end{cases}
     $$
  Put
  $$ f_t := p_1^{u_1/2}p_2^{u_2/2}...p_n^{u_n/2}.$$
  Then $q = q_1^2\in \mathcal{S}(Q_0)$ is divisible by $t$ iff $q_1$ is divisible by $f_t$. Thus, writing $q_2:=q_1/f_t$, we have 
  \begin{equation*}
  \mathcal{S}_t(Q_0) = \left\{q_2^2g_t : tg_tq_2^2\in\mathcal{S}(Q_0)\right\},
  \end{equation*}
  where 
  $$
  g_t := \frac{f_t^2}{t} = \epsilon^{-1}p_1^{u_1-v_1}...p_n^{u_n-v_n}.
  $$
  We observe that 
  \begin{equation*}
  \begin{split}
  tg_tq_2^2\in \mathcal{S}(Q_0) & \Longleftrightarrow Q_0<\mathcal{N}(tg_tq_2^2)\le 2Q_0 \\ &
                        \Longleftrightarrow Q_0< \mathcal{N}(f_t^2q_2^2)\le 2Q_0 \\& 
                        \Longleftrightarrow \frac{\sqrt{Q_0}}{\mathcal{N}(f_t)}< \mathcal{N}(q_2)\le \frac{\sqrt{2Q_0}}{\mathcal{N}(f_t)}.
 \end{split}
  \end{equation*}
  Hence,
  $$\mathcal{S}_t(Q_0) = \left\{q_2^2g_t : \frac{\sqrt{Q_0}}{\mathcal{N}(f_t)}< \mathcal{N}(q_2)\le \frac{\sqrt{2Q_0}}{\mathcal{N}(f_t)} \right\}.  
  $$
  
   As previously, we suppose that $0\le u\leq \sqrt{2Q_0}/t$, $k\in \mathbb{Z}[i]\setminus\{0\}$, $l\in \mathbb{Z}[i]$ and $(k,l) =1$ and define 
   \begin{equation*}
   \begin{split}
  A_t(u,k,l) := \sup\limits_{|y|\le \frac{\sqrt{2Q_0}}{|t|}}
  \left|\left\{q\in\mathcal{S}_t(Q_0)\cap B(y,u) : q\equiv l\bmod{k}\right\}\right|. 
   \end{split}
   \end{equation*}
   Noting that $q> \sqrt{Q_0}/|t|$ if $q\in \mathcal{S}_t(Q_0)$, it follows that
   $$
   A_t(u,k,l)\le  \sup\limits_{\frac{\sqrt{Q_0}}{|t|}\le |y|\le \frac{\sqrt{2Q_0}}{|t|}} |\mathcal{A}(y)|,
   $$
   where 
   $$
   \mathcal{A}(y):= \left\{ 
                 q_2\in B\left(0,\left(\frac{\sqrt{2Q_0}}{|tg_t|}\right)^{1/2}\right)\ :\ q_2^2g_t\in B(y,u), \ 
                 q_2^2g_t\equiv l\bmod{k}\right\}.
                 $$
 
 We aim to verify the condition \eqref{condi} for $X=N^{\varepsilon}$. Thus, our next task is to bound the cardinality of $\mathcal{A}(y)$ if
 $$
 \frac{\sqrt{Q_0}}{|t|}\le |y|\le \frac{\sqrt{2Q_0}}{|t|}.
 $$
 Let $\delta_t(k,l)$ be the number of solutions $x\bmod{k}$ to the congruence
 $$ x^2g_t\equiv l\bmod{k}.$$
 Then the number of Gaussian integers $x$ contained in a ball $B(a,r)$ and satisfying the congruence $x^2g_t\equiv l\bmod{k}$ is 
 $$\ll \left(1+\frac{r^2}{\mathcal{N}(k)}\right)\delta_t(k,l).$$
 We deduce that 
 \begin{equation} \label{theb}
 \mathcal{A}(y)\ll \left(1+ \frac{\sqrt{Q_0}}{|tg_t|\mathcal{N}(k)}\right)\delta_t(k,l)
 \end{equation}
 by ignoring the condition that $q_2^2g_t\in B(y,u)$.
 We shall use \eqref{theb} if $\sqrt{2Q_0}/(2|t|)<u\le \sqrt{2Q_0}/|t|$. If $u\le \sqrt{2Q_0}/(2|t|)$, then we 
 obtain a stronger bound as follows. Note that
 $$
 q_2^2g_t\in B(y,u) \Longleftrightarrow \left| q_2^2-\frac{y}{g_t} \right|\le \frac{u}{|g_t|}.
 $$
 Hence,
 \begin{equation} \label{spli}
 \begin{split}
 \left|q_2-\sqrt{\frac{y}{g_t}}\right|\cdot \left|q_2+\sqrt{\frac{y}{g_t}}\right|\le \frac{u}{|g_t|},
 \end{split}
 \end{equation}
 where we define the square root of the complex number $s=\rho e^{i\phi}$ to be $\rho^{1/2}e^{i\phi/2}$ if $\rho:=|s|$ and $-\pi<\phi=\arg(s)\le \pi$.   
 The set 
 $$
 \left\{q_2 : \left|q_2^2 - \frac{y}{g_t}\right|\le \frac{u}{|g_t|}\right\}
 $$
 consists of two connected components, one containing $\sqrt{y/g_t}$ and the other one containing $-\sqrt{y/g_t}$. By symmetry, it suffices to look at the case when $q_2$ is contained 
 in the first component. Moreover, we may restrict ourselves to the case when $0\le \arg(y/g_t)\le \pi/2$ since all other cases are similar. In this case, $\Re(q_2)\ge 0$ and hence,
 \begin{equation} \label{upp}
 \begin{split}
 \left|q_2+\sqrt{\frac{y}{g_t}}\right|\geq \Re\left(q_2+\sqrt{\frac{y}{g_t}}\right) 
 \geq \Re \left(\sqrt{\frac{y}{g_t}}\right) \geq \sqrt{\left|\frac{y}{2g_t}\right|}.
 \end{split}
 \end{equation}
Combining \eqref{spli} and \eqref{upp}, we get
\begin{equation*}
\begin{split}
\left|q_2-\sqrt{\frac{y}{g_t}}\right| \le \frac{\sqrt{2}u}{\sqrt{|yg_t|}} 
\end{split}
\end{equation*}
and hence 
$$
q_2\in B\left(\sqrt{\frac{y}{g_t}},\frac{\sqrt{2}u}{\sqrt{|yg_t|}}\right).
$$
Thus, 
\begin{equation*}
\mathcal{A}(y) \ll \left(1+\frac{u^2}{|yg_t|\mathcal{N}(k)}\right)\delta_t(k,l),
\end{equation*}
which contains also the bound \eqref{theb} for the case when $\sqrt{2Q_0}/(2|t|)<u\le \sqrt{2Q_0}/|t|$. We observe that
\begin{equation*}
\begin{split}
             \frac{|\mathcal{S}_t(Q_0)|/\mathcal{N}(k)}{Q_0/|t|^2}u^2 \ge c_6\frac{u^2}{|yg_t|\mathcal{N}(k)}.
\end{split}
\end{equation*}
Combining everything, it follows that
\begin{equation*}
A_t(u,k,l)\le c_{7}\left(1+\frac{|\mathcal{S}_t(Q_0)|/\mathcal{N}(k)}{Q_0/|t|^2}u^2\right)\delta_t(k,l),
\end{equation*}
which is \eqref{condi} with $X$ replaced by $\delta_t(k,l)$.

The remaining task is to bound $\delta_t(k,l)$.
If $(g_t,k)>1$, then $\delta_t(k,l) =0$ since $k$ and $l$ are supposed to be coprime. 
Therefore, we can assume that $(g_t,k) =1$. Let $g\bmod{k}$ be a multiplicative inverse of $g_t\bmod{k}$,
i.e. $gg_t\equiv 1\bmod{k}$. Put $l^*= gl$. Then $(k,l^*)=1$, and $x^2g_t\equiv l\bmod{k}$ is equivalent to $x^2\equiv l^*\bmod{k}$.
Now using the Chinese remainder theorem, $\delta_t(k,l)$ equals the product of all $\delta_t(p^{a(p)},l)$, where 
$$
k=\epsilon \prod\limits_{p|k} p^{a(p)}
$$
is a prime factorization of $k$ in $\mathbb{Z}[i]$. Thus it suffices to look at prime powers $k=p^a$. If $a=1$ and $p\nmid 2$, then since our congruence is quadratic, we deduce that 
$\delta_t(p^a,l)\le 2$ by looking at the residue field of $(p)$. If $a>1$ and $p\nmid  2$, then using Hensel's Lemma for $\mathbb{Z}[i]$, this inequality remains valid. 
If $p|2$, an application of Hensel's lemma shows that $\delta_t(2^a,l)\le 4$.  It follows that for all $k\in\mathbb{Z}[i]$, 
 $$ \delta_t(k,l)\le 2^{\omega(k)+1},$$
 where $\omega(k)$ is the number of distinct prime divisors of $k$. Furthermore, since $\mathcal{N}(k)\le \sqrt{N}$, we have
 $$ 2^{\omega(k)} \le d(\mathcal{N}(k)) \ll N^{\varepsilon},$$
 where $d(n)$ denotes the number of positive divisors of the integer $n$.
 Therefore, \eqref{condi} holds with 
 \begin{equation}
  X := c_{7}N^{\varepsilon}.
  \end{equation}
Now combining Theorem \ref{result2}, Lemma \ref{cr3} and $|\mathcal{S}|\ll \sqrt{Q_0}$, we obtain
$$
\sum\limits_{q\in\mathcal{S}(Q_0)}\sum\limits_{\substack{a\bmod{q}\\ (a,q)=1}}\left|S\left(\frac{a}{q}\right)\right|^2\le c_{8}\left(N+Q_0N^\varepsilon\left(\sqrt{N}+\sqrt{Q_0}\right)\right)Z.
$$
This implies the desired result by adding the contributions of dyadic intervals. $\Box$

\section{Proof of Theorem \ref{primemoduli}}
 We prove Theorem \ref{primemoduli} using Theorem \ref{result1} with $\mathcal{S}$ being the set of all Gaussian primes with norm $\le Q$. 
 Thus we need to estimate the term $A_t(u,k,l)$. First we consider the case when $t$ is a unit in $\mathbb{Z}[i]$. 
 By the Brun-Titchmarsh inequality for number fields developed in \cite{Sch}, we have 
 \begin{equation*}
 A_t(u,k,l) \le \frac{2^5}{\pi}\frac{u^2}{\varphi(k)\log\left(u^2/\mathcal{N}(k)\right)}
 \end{equation*}
 if $u^2/\mathcal{N}(k) > 1$. It follows that 
 \begin{equation} \label{BrunTitchmarsh}
 A_t\left(\frac{\sqrt{Q}}{N^{1/2}|z|},\frac{r}{t},l\right) \le \frac{2^5}{\pi}\frac{Q}{|z|^2 N\varphi(r)\log\left(Q/({\mathcal{N}(r)}|z|^2N)\right)}
 \end{equation}
 if $Q/(N|z|^2) > \mathcal{N}(r)$.  
 Here, $\varphi(n)$ is the Euler function for Gaussian integers, defined to be the number of residue
 classes $a \bmod{n}$ in $\mathbb{Z}[i]$ such that $(a,n)=1$.
 If $|rz|\le 2N^{-1/4}$, then $Q/(N|z|^2) > \mathcal{N}(r)$ is satisfied since $4/\sqrt{N}< Q/N$ by the condition given in Theorem \ref{primemoduli}. 
 From the assumptions of Theorem \ref{primemoduli} and \eqref{BrunTitchmarsh}, we deduce that
 \begin{equation}
 \sum\limits_{\substack{0< |m|\le 3|rz|\sqrt{Q}\\ (m,\frac{r}{t}) = 1}} 
 A_t\left(\frac{\sqrt{Q}}{|z|N^{1/2}},\frac{r}{t},hm\right)\\ \le \frac{9\cdot 2^5}{\pi}\frac{2Q^2}{(1-\delta)N\log Q}
 \cdot\frac{\mathcal{N}(r)}{\varphi{(r)}}
 \end{equation} 
 for any $h\in \mathbb{Z}[i]$ with $(h,r) = 1$. It is well-known that
 \begin{equation}
 \begin{split}
 \frac{\mathcal{N}(r)}{\varphi{(r)}} \ll \log\log\mathcal{N}(r).
 \end{split} 
 \end{equation}
 Hence, 
 \begin{equation} \label{last1}
 \begin{split}
 \sum\limits_{\substack{0< |m|\le 3|rz|\sqrt{Q}\\ (m,\frac{r}{t}) = 1}} A_t\left(\frac{\sqrt{Q}}{|z|N^{1/2}},\frac{r}{t},hm\right) 
 \le c_9 \frac{Q^2\log\log\mathcal{N}(r)}{N(1-\delta)\log Q}.
 \end{split}
 \end{equation}
 
If $t$ is non-zero and a non-unit Gaussian integer, then $\mathcal{S}_t$ contains at most $4$ elements. Hence, the contribution of these $t$'s is bounded by 
\begin{equation} \label{last2}
\begin{split}
 \sum\limits_{\substack{t|r\\ \mathcal{N}(t) \geq 2}}\sum\limits_{\substack{0<|m|\le \frac{3|rz|\sqrt{Q}}{|t|} \\ (m,\frac{r}{t}) = 1}} 
 A_t\left(\frac{\sqrt{Q}}{\sqrt{N}|z||t|},\frac{r}{t}, hm \right)
 \le & \sum\limits_{t|r} \frac{36|rz|^2Q}{|t|^2}\\  \le & c_{10}|rz|^2Q \log\log \mathcal{N}(r).
 \end{split}
\end{equation}
 Combining Theorem \ref{result1}, \eqref{last1} and \eqref{last2}, 
we obtain the claimed result by a short calculation upon noting that 
$$
\frac{Q^{(1-\delta)/2}}{1-\delta} \ge \frac{1}{2}\log Q. \ \Box
$$

 \end{document}